\documentclass[secthm,seceqn,amsthm,ussrhead,10pt]{amsart}
\usepackage{amsmath,latexsym}
\usepackage[english]{babel}
\usepackage[psamsfonts]{amssymb}
\usepackage{times}
\usepackage{cite}

\usepackage[mathcal]{euscript}

\newtheorem{Th}{Theorem}
\newtheorem{Lem}[Th]{Lemma}
\newtheorem{proposition}[Th]{Lemma}

\begin{document}
\sloppy
\thispagestyle{empty}

\textbf{The structure of simple noncommutative Jordan superalgebras}
\footnote{
The results from the third section of the paper (Theorems 5-8) were obtained by the second author with financial support by RFBR No. 16-31-60111 (mol\underline{ }a\underline{ }dk). The results from the 4 and 5 sections (Theorems 9-12) were obtained by the first and the third authors with financial support by CNPq No 300603/2016-9; RFBR 17-01-00258 and the President’s Programme ‘Support of Young Russian Scientists’ (grant MK-1378.2017.1).}

\medskip

\medskip
\textbf{Ivan Kaygorodov$^{a}$, Artem Lopatin$^{b}$, Yury Popov$^{c}$}

\medskip
{\tiny

$^a$ Universidade Federal do ABC, CMCC, Santo Andre, SP, Brazil

$^b$ Sobolev Institute of Mathematics, Omsk Branch, SB RAS, Omsk, Russia

$^c$ Universidade Estadual de Campinas, IMECC, Campinas, SP, Brazil


    E-mail addresses:\smallskip

    Ivan Kaygorodov (kaygorodov.ivan@gmail.com)
    
    Artem Lopatin (dr.artem.lopatin@gmail.com)

    Yury Popov (yuri.ppv@gmail.com)

}    

\section*{Abstract}

In this paper we describe all subalgebras and automorphisms of simple noncommutative Jordan superalgebras $K_3(\alpha,\beta,\gamma)$ and $D_t(\alpha,\beta,\gamma);$ 
and compute the derivations of the nontrivial simple finite-dimensional noncommutative Jordan superalgebras.

\medskip

\section{Introduction}
The class of noncommutative Jordan algebras is vast: 
for example, it includes alternative algebras, Jordan algebras, quasiassociative algebras, quadratic flexible algebras and anticommutative algebras. 
Schafer proved that a simple noncommutative Jordan algebra is either a simple Jordan algebra, or a simple quasiassociative algebra, or a simple flexible algebra of degree 2 \cite{Sch}. 
Oehmke proved the analog of Schafer's classification for simple flexible algebras with strictly associative powers of and of characteristic $\neq 2, 3$ \cite{Oeh}, 
McCrimmon classified simple noncommutative Jordan algebras of degree $>2$ and characteristic $\neq 2$ \cite{McC,mcC2}, and  
Smith described such algebras of degree 2 \cite{Smith}. 
The case of nodal simple algebras of positive characteristic was considered in the papers of Kokoris \cite{Kok,Kok2}.

The case of simple finite-dimensional Jordan superalgebras over algebraically closed fields of characteristic $0$ was studied by Kac \cite{Kac1} and Kantor \cite{Kan}. 
Racine and Zelmanov \cite{RZ} classified the finite-dimensional Jordan superalgebras of characteristic $\neq 2$ with semisimple even part, the case when even part is not semisimple was considered by Martinez and Zelmanov in \cite{MZ01}
and Cantarini and Kac described all linearly compact simple Jordan superalgebras \cite{KacKant07}. 
Simple noncommutative Jordan superalgebras were described by Pozhidaev and Shestakov in \cite{ps1,ps2}.
Representations of simple noncommutative superalgebras were described by Yu. Popov \cite{Po1}.

Nowadays, a great interest is shown in the study of nonassociative algebras and superalgebras with derivations;
in particular, Jordan algebras and superalgebras.
For example, A. Popov determined the structure of differentiably simple Jordan algebras \cite{Popov2};
Kaygorodov and Yu. Popov determined the structure of Jordan algebras admitting invertible Leibniz derivations \cite{KP,KP16};
Kaygorodov, Lopatin and Yu. Popov determined the structure of Jordan algebras admitting derivations with invertible values \cite{KP2};
Barreiro, Elduque and Mart\'{i}nez described the derivations of Cheng-Kac Jordan superalgebra \cite{BEC}; 
Kaygorodov,  Shestakov and  Zhelyabin studied generalized derivations of Jordan algebras and superalgebras  \cite{kay07,kay10,kay11,kay12,shestakov2}.
Another interesting problem in the study of Jordan algebras and superalgebras is a description of maximal subalgebras and automorphisms \cite{ELS,ELS2}.

In this paper we are studying simple noncommutative Jordan superalgebras constructed in some papers of Pozhidaev and Shestakov.
We describe all subalgebras and automorphisms of simple noncommutative Jordan superalgebras $K_3(\alpha,\beta,\gamma)$ and $D_t(\alpha,\beta,\gamma);$ 
and compute the derivations of the nontrivial simple finite-dimensional noncommutative Jordan superalgebras.

\medskip

\medskip

\section{Preliminaries}
Let $U = U_{\bar{0}} + U_{\bar{1}}$ be a superalgebra, $(-1)^{xy} = (-1)^{p(x)p(y)},$ where $p(x)$ is the parity of $x,$ ($p(x) = i,$ if $x \in U_{\bar{i}}$).
In what follows, if the parity of an element arises in a formula, this element is assumed to be homogeneous.
By $L_x, R_x$ we denote the operators of left and right multiplication by $x \in U:$
$$yL_x = (-1)^{xy}xy, yR_x = yx; [x,y] = xy - (-1)^{xy}yx, x \bullet y = xy + (-1)^{xy}yx.$$

A commutative superalgebra $J$ is called \emph{Jordan superalgebra} if it satifies the following identity:
\begin{equation}
\label{jord_identity}
R_aR_bR_c + (-1)^{a,b,c}R_cR_bR_a + (-1)^{bc}R_{(ac)b} = R_aR_{bc} + (-1)^{a,b,c}R_cR_{ab} + (-1)^{ab}R_bR_{ac}.
\end{equation}

A superalgebra $U$ is called a \emph{noncommutative Jordan superalgebra} if it satisfies the following operator identities:
\begin{equation}
\label{noncomm_jord_identity}
[R_{x \circ y}, L_z] + (-1)^{x(y+z)}[R_{y\circ z}, L_x] + (-1)^{z(x+y)}[R_{z \circ x}, L_y] \end{equation}
\begin{equation}
[R_x, L_y] = [L_x, R_y].\label{flex}\end{equation}
The second operator identity defines the class of \emph{flexible superalgebras.}
If we assume that all elements are even we arrive at the notion of a noncommutative Jordan algebra.


A homogeneous linear mapping $d$ of $U$ is called a \textit{derivation} of $U$, if it satisfies the following relation:
\begin{equation}
\label{der_con}
(xy)d = (-1)^{dy} xd\cdot y + x \cdot yd.
\end{equation}
\medskip

A binary linear operation $\{,\}$ is called a \emph{generic Poisson bracket} \cite{ksu} on a superalgebra $(A,\cdot)$ if for arbitrary homogeneous $a, b, c \in A$ we have
\begin{equation}
\label{pois_br}
\{a \cdot b, c\} = (-1)^{bc}\{a,c\}\cdot b + a \cdot \{b,c\}.
\end{equation}

We notice that there is a one-to-one correspondence between noncommutative Jordan superalgebras and superanticommutative Poisson brackets on adjoint Jordan superalgebras:
\begin{Lem}
\label{JordPois}
\cite{ps2}
Let $(J, \bullet)$ be a Jordan superalgebra and $[,]$ be a superanticommutative Poisson bracket on $J.$ Then the operation $ab = \frac{1}{2}(a\bullet b + [a,b])$ turns $J$ into a noncommutative Jordan superalgebra. Conversely, if $U$ is a noncommutative Jordan superalgebra, then the supercommutator $[,]$ is a Poisson bracket on a Jordan superalgebra $U^{(+)}.$ Moreover, the multiplication in $U$ can be recovered by the Jordan multiplication in $U^{(+)}$ and the Poisson bracket $[,]: ab = \frac{1}{2}(a \bullet b + [a,b]).$
\end{Lem}

We provide examples of noncommutative Jordan superalgebras given in \cite{ps1,ps2}:

\subsection{ The superalgebra $K_3(\alpha, \beta, \gamma)$}
The superalgebra 
$K_3(\alpha, \beta, \gamma) = U_{\bar{0}} \oplus U_{\bar{1}}, U_{\bar{0}} = \langle e \rangle, U_{\bar{1}} = \langle z, w \rangle$ 
is defined by the following multiplication table:
\begin{table}[h]
\begin{tabular}{c|c|c|c}
    & $e$ & $z$ & $w$ \\
    \hline
$e$&         $e$                     & $\alpha z + \beta w$ & $\gamma z + (1-\alpha)w$ \\ 
\hline
$z$& $(1-\alpha)z - \beta w$ & $-2 \beta e$         & $2 \alpha e$             \\ 
\hline
$w$& $\alpha w - \gamma z$   & $-2(1- \alpha)e$     & $2 \gamma e$             \\ 
\end{tabular}
\end{table}

The superalgebra $K_3(\alpha, \beta, \gamma)^{(+)}$ is isomorphic to the simple nonunital Jordan superalgebra $K_3 = K_3(\frac{1}{2}, 0 , 0).$
In \cite{Po1} was proved the following result:

\begin{Lem}
\label{K3_class}
\textit{If $\mathbb{F}$ is a field which allows square root extraction, then for $\alpha, \beta, \gamma \in \mathbb{F},$  
$K_3(\alpha, \beta, \gamma)$ is isomorphic either to $K_3(\lambda,0,0) = K_3(\lambda)$ for some $\lambda \in \mathbb{F}$, or to $K_3(\frac{1}{2},\frac{1}{2},0)=K_3^{1/2}.$}
\end{Lem}


\subsection{The superalgebra $D_t(\alpha, \beta, \gamma)$}
Take $\alpha, \beta, \gamma, t \in \mathbb{F}.$ We define the superalgebra 
$U = D_t(\alpha, \beta, \gamma)=U_{\bar{0}} \oplus U_{\bar{1}},$ 
$U_{\bar{0}} = \langle e_1, e_2 \rangle , U_{\bar{1}} = \langle x, y \rangle$
is defined by the following multiplication table:

$$\begin{array} {c|c|c|c|c}
    & e_1 & e_2 & x & y \\ \hline
e_1 & e_1                   & 0                           & \alpha x + \beta y    &  \gamma x + (1-\alpha)y \\ \hline
e_2 & 0                     & e_2                         & (1-\alpha)x - \beta y & -\gamma x + \alpha y \\ \hline
x   & (1-\alpha)x - \beta y & \alpha x + \beta y          & -2\beta(e_1 - te_2)   & 2(\alpha e_1 + (1-\alpha)te_2) \\ \hline
y   & -\gamma x + \alpha y  & \gamma x + (1-\alpha)y      &  -2((1-\alpha)e_1 + \alpha t e_2) & 2\gamma(e_1 - te_2) \\
\end{array}$$

Putting $t = -1,$ we obtain the superlagebra $M_{1,1}(\alpha, \beta, \gamma)$, and putting $t = -2$, we obtain the superalgebra $osp(1,2)(\alpha, \beta, \gamma)$ (see \cite{ps2}).
The superalgebra $D_t(\alpha, \beta, \gamma)^{(+)}$ is isomorphic to the Jordan superalgebra $D_t = D_t(\frac{1}{2})$.
One can easily see that nonsimple noncommutative Jordan superalgebra $D_0(\alpha, \beta, \gamma)$ is a unital hull of $K_3(\alpha, \beta, \gamma)$.

In \cite{Po1} was proved the following result:

\begin{Lem}
\label{Dt_class}
\textit{If $\mathbb{F}$ is a field which allows square root extraction, then for $\alpha, \beta, \gamma \in \mathbb{F},$  $D_t(\alpha, \beta, \gamma)$ is isomorphic either to $D_t(\lambda,0,0) = D_t(\lambda)$ for some $\lambda \in \mathbb{F}$, or to $D_t(\frac{1}{2},\frac{1}{2},0)=D_t^{1/2}.$}
\end{Lem}

\subsection{The superalgebra $U(V, f, \star)$} 
Let $V = V_0 \oplus V_1$ be a vector superspace over  $\mathbb{F},$ and let $f$ be a supersymmetric nondegenerate bilinear form on $V.$ Also let $\star$ be a superanticommutative multiplication on $V$ such that $f(x \star y, z) = f(x, y \star z).$ Then we can define a multiplication on $U = \mathbb{F} \oplus V$ in the following way:
$$(\alpha + x)(\beta + y) = (\alpha \beta + f(x,y)) + (\alpha y + \beta x + x \star y),$$
and the resulting superalgebra is denoted $U(V, f, \star).$\\

The superalgebra $U(V, f, \star)^{(+)}$ is isomorphic to the simple Jordan superalgebra $ J(V,f)$ of nondegenerate supersymmetric bilinear form. Also note that $J(V,f) = U(V,f,0).$

\subsection{The superalgebra $J(\Gamma_n, A)$}
Let $\Gamma$ be the Grassmann superalgebra in generators $1, x_i, i \in \mathcal{I},$ where $\mathcal{I}$ can be empty, and $\Gamma_n$ be the Grassmann superalgebra in generators $1, x_1, \dots, x_n.$

We define the new operation $\{,\}$ on $\Gamma$ (\emph{the Poisson-Grassmann bracket})  by defining for $f, g \in \Gamma_{\bar{0}} \cup \Gamma_{\bar{1}}$
\begin{equation}
\label{Pois_br}
\{f,g\} = (-1)^f \sum_{j=1}^{\infty} \frac{\partial f}{\partial x_j} \frac{\partial g}{\partial x_j},
\end{equation}
where
\begin{equation*}
\frac{\partial}{\partial x_j} (x_{i_1} x_{i_2} \dots x_{i_n}) = \begin{cases} (-1)^{k-1}x_{i_1}x_{i_2} \dots x_{i_{k-1}}x_{i_{k+1}}\dots x_{i_n}, & \mbox{if } j = i_k,\\ 0, & \mbox{if } j \notin {i_1, i_2, \dots, i_n}. \end{cases}
\end{equation*}
Let ${\bar{\Gamma}}$ be the isomorphic copy of $\Gamma$ with the isomorphism mapping $x \rightarrow \bar{x}.$ We define the structure of a Jordan superalgebra  on $J(\Gamma) = \Gamma \oplus \bar{\Gamma}$, by setting $J(\Gamma)_{\bar{0}} = \Gamma_{\bar{0}} + \bar{\Gamma}_{\bar{1}}, J(\Gamma)_{\bar{1}} = \Gamma_{\bar{1}} + \bar{\Gamma}_{\bar{0}}$ and defining multiplication by the rule
$$a\cdot b = ab, \bar{a}\cdot b = (-1)^b\overline{ab}, a\cdot \bar{b} = \overline{ab}, \bar{a} \cdot \bar{b} = (-1)^b \{a,b\},$$
where $a, b \in \Gamma_{\bar{0}} \cup \Gamma_{\bar{1}}$ and $ab$ is their product in $\Gamma$. By $J(\Gamma_n)$ we will denote the subalgebra $\Gamma_n + \bar{\Gamma}_n$ of $J(\Gamma)$.

If $\Gamma_n$ is considered as a Poisson superalgebra, then it is easily seen that $J(\Gamma_n)$ is the \emph{Kantor double} of $\Gamma_n.$ Kantor doubles and their derivations were studied in many works, such as \cite{Kan, Ret, kay11}, and many others.
Let $A \in (\Gamma_n)_{\bar{0}}.$ Define a superanticommutative binary bilinear operation $\{, \}$ on $J(\Gamma_n)$ by the rule
$$\{\bar{a}, \bar{b}\} = (-1)^b abA,$$
and zero otherwise. One can check that this operation is a Poisson bracket on $J(\Gamma_n).$ On the underlying vector space of $J(\Gamma_n)$ we define a binary operation
$$ab = a\cdot b + \{a,b\}$$
The resulting superalgebra will be denoted by $J(\Gamma_n, A).$

It is easy to see that $J(\Gamma_n,A)^{(+)} = J(\Gamma_n) = J(\Gamma_n,0).$

\subsection{The superalgebra $\Gamma_n(\mathcal{D})$}
Let $\mathcal{D} = \{d_i, i = 1, \dots, n\}$ be a set of odd derivations of Grassmann superalgebra $(\Gamma_n, \cdot)$ satisfying the condition $x_i d_j = x_j d_i = a_{ij} \in (\Gamma_n)_{\bar{0}}, i, j = 1,\dots, n$, such that $\Gamma_n$ is differentiably simple with respect to $\mathcal{D}.$ We define a Poisson bracket on $\Gamma_n$ by the rule $[x_i,x_j] = x_i d_j$.

On the vector space of $\Gamma_n$ we define a new multiplication by $ab = a\cdot b + [a,b]$. The resulting noncommutative Jordan superalgebra is denoted by $\Gamma_n(\mathcal{D}).$

By construction, the superalgebra $\Gamma_n(\mathcal{D})^{(+)}$ is isomorphic to the Grassmann superalgebra $\Gamma_n.$
\medskip

Pozhidaev and Shestakov described all simple central noncommutative Jordan superalgebras:
\begin{Th}\cite{ps2}
\label{classification}
Let $U$ be a simple noncommutative central Jordan superalgebra over a field $\mathbb{F}$ of characteristic zero. 
Suppose that $U$ is neither supercommutative nor quasiassociative. Then one of the following holds$:$

\begin{enumerate}
    
\item[$\bullet$] $U \cong K_3(\alpha, \beta, \gamma), \ M_{1,1}(F)^q, \  osp (1,2)^q, \ J(\Gamma_n, A), \ \Gamma_n(\mathcal{D});$
\end{enumerate}

or there exists an extension $P$ of $\mathbb{F}$ of degree $\leq 2$ such that $U \otimes_{\mathbb{F}} P$ is isomorphic as a $P$-superalgebra to one of the following ones$:$
\begin{enumerate}

\item[$\bullet$] $ D_t(\alpha), \ U(V, f, \star).$
\end{enumerate}
\end{Th}

\medskip

\section{Subalgebras of $K_3(\alpha,\beta,\gamma)$ and $D_t(\alpha,\beta,\gamma)$}

\subsection{Subalgebras of $K_3(\alpha)$}
The study of subalgebras of the superalgebra $K_3(\alpha)$ has $2$ parts.
We are considering $1$-dimensional and $2$-dimensional subalgebras.
The result of this subsection is the following 

\begin{Th} \label{subk}
Let $M$ be a subalgebra of $K_3(\alpha),$ then
\begin{enumerate}

\item if $dim(M)=1$ and $\alpha=\frac{1}{2}$ then $M$ is generated by one of the following vectors
$\gamma_1e+\gamma_2z+\gamma_3w,$ where $\gamma_1,\gamma_2,\gamma_3 \in \mathbb{F};$

\item if $dim(M)=1$ and $\alpha\neq\frac{1}{2}$ then $M$ is generated by one of the following vectors
$\gamma_1e+\gamma_2z$ or $\gamma_1e+\gamma_2w,$ where $\gamma_1,\gamma_2 \in \mathbb{F};$

\item if $dim(M)=2$ and $\alpha=\frac{1}{2}$ then $M$ is generated by one of the following pairs of vectors
$e$ and $\gamma_1 z+\gamma_2 w,$ where $\gamma_1,\gamma_2 \in \mathbb{F};$

\item if $dim(M)=2$ and $\alpha\neq\frac{1}{2}$ then $M$ is generated by one of the following pairs of vectors
$e$ and $z$ or $e$ and $w.$

\end{enumerate}

\end{Th}

The proof of the Theorem \ref{subk} is based on the following subsections \ref{k1} and \ref{k2}.

\subsubsection{$1$-dimensional subalgebras}\label{k1}
Let $M$ be an $1$-dimensional subalgebra of $K_3(\alpha),$ then 
$M$ is generated by the following vector
$w_1=x_1 e+x_2 z+x_3 w,$ for some $x_1,x_2,x_3 \in \mathbb{F}.$
Now, 
$$w_1^2=x_1^2 e+x_1x_2z+x_1x_3w+(4\alpha-2)x_2x_3 e$$ and for an element $k \in \mathbb{F}$ we have
$$x_1x_2=kx_2, \ x_1x_3=kx_3, \ x_1^2+(4\alpha-2)x_2x_3=kx_1.$$
It is easy to see, that if $\alpha=\frac{1}{2},$ then $w_1^2=kw_1$ and for any vector $w_1$ we have $1$-dimensional subalgebra.

Now, if $\alpha\neq\frac{1}{2},$ then
\begin{enumerate}
\item[$\bullet$] if $x_2\neq0,$ then $x_1=k, x_3=0$ and $M$ is generated by some element $x_1e+x_2z;$

\item[$\bullet$] if $x_2=0,$ then $x_1^2=kx_1$ and $M$ is generated by some element $x_1e+x_3w.$

\end{enumerate}

\subsubsection{$2$-dimensional subalgebras}\label{k2}
Let $M$ be a $2$-dimensional subalgebra of $K_3(\alpha),$ then 
$M$ is generated by the following vectors
$w_1=x_1 e+x_2 z+x_3 w$ and $w_2=y_1 e+y_2 z+y_3 w.$ 
It is easy to see, that we can suppose 
$w_1=e+x_2 z+x_3 w$ and $w_2=y_2 z+y_3 w.$ 
Here, we have two separate cases.
\begin{enumerate}
    \item[I.]If $x_2 z+x_3 w$ and $y_2 z+y_3 w$ are linearly dependent, we have that $M$ is generated by $e$ and $x_2 z+x_3 w.$
Now, $ew_2=\alpha y_2 z+(1-\alpha)y_3 w \in M$ and if $y_2\neq0, y_3\neq0$ and $\alpha\neq\frac{1}{2},$ then $dim(M)=3.$

On the other hand,

\begin{enumerate}

\item[$\bullet$] if $\alpha=\frac{1}{2},$ then $e$ and $x_2 z+x_3 w$ are generating a $2$-dimensional subalgebra.

\item[$\bullet$] if $y_2=0,$ then $M$ is generated by $e$ and $w;$

\item[$\bullet$] if $y_3=0,$ then $M$ is generated by $e$ and $z.$
\end{enumerate}

\item[II.] If $x_2 z+x_3 w$ and $y_2 z+y_3 w$ are linearly independent, we have
$w_2^2=y_2y_3(4\alpha-2)e \in M.$

If $\alpha=\frac{1}{2},$ then $w_1w_2=\frac{1}{2}w_2+(x_2y_2-x_3y_3)e \in M$ and $e\in M.$ It is follows, that $dim(M)=3.$

If $\alpha\neq \frac{1}{2},$ then we have three opportunities: $y_2=0, y_3=0$ or $e\in M.$

\begin{enumerate}

\item[$\bullet$]If $y_2=0,$ we have that $w \in M$ and $e+x_2z\in M.$
It is $$(e+x_2z) w+w(e+x_2z)= w + (4\alpha-2) x_2 e \in M$$ and $dim(M)=3.$

\item[$\bullet$]If $y_3=0$ or $e\in M$ by some similar way we have that $dim(M)=3.$
\end{enumerate}

\end{enumerate}

\subsection{Subalgebras of $K_3^{1/2}$}
The study of subalgebras of the superalgebra $K_3^{1/2}$ has $2$ parts.
We are considering $1$-dimensional and $2$-dimensional subalgebras.
The result of this subsection is the following 

\begin{Th}\label{subk12}
Let $M$ be a subalgebra of $K_3^{1/2},$ then
\begin{enumerate}

\item if $dim(M)=1,$ then $M$ is generated by the following vector $\gamma_1e+\gamma_2w,$ where  $\gamma_1,\gamma_2 \in \mathbb{F};$

\item if $dim(M)=2,$ then $M$ is generated by the following  vectors $e$ and $w.$

\end{enumerate}

\end{Th}

The proof of the Theorem \ref{subk12} is based on the following subsections \ref{k121} and \ref{k122}.

\subsubsection{$1$-dimensional subalgebras}\label{k121}
Let $M$ be an $1$-dimensional subalgebra of $K_3^{1/2},$ then 
$M$ is generated by the following vector
$w_1=x_1 e+x_2 z+x_3 w.$
Now, 
$w_1^2=x_1^2 e+x_1x_2z+x_1x_3w- x_2^2e$ and for an element $k \in \mathbb{F}$ we have
$$kx_1=x_1^2-x_2^2, \ kx_2=x_1x_2, \ kx_3=x_1x_3.$$

\begin{enumerate}
\item[$\bullet$] If $x_2\neq0,$ then $x_1=k$ and $x_2=0.$
                
\item[$\bullet$] If $x_2=0,$ then $x_1^2=kx_1$ and $M$ is generated by  $x_1e+x_3w.$
\end{enumerate}

\subsubsection{$2$-dimensional subalgebras}\label{k122}
Let $M$ be a $2$-dimensional subalgebra of $K_3^{1/2},$ then 
$M$ is generated by the following vectors
$w_1=x_1 e+x_2 z+x_3 w$ and $w_2=y_1 e+y_2 z+y_3 w.$ 
It is easy to see, that we can suppose that 
$w_1=e+x_2 z+x_3 w$ and $w_2=y_2 z+y_3 w.$ 

\begin{enumerate}
    \item[I.]If $x_2 z+x_3 w$ and $y_2 z+y_3 w$ are linearly dependent, we have that 
$e \in M$ and $2ew_2= w_2+ y_2w \in M.$
Now, $M$ is generated by $e$ and $w.$

\item[II.] If  $x_2 z+x_3 w$ and $y_2 z+y_3 w$ are linearly independent, 
we have that $w_2^2=-y_2^2e \in M.$
Here, it is true that $e\in M$ or $y_2=0,$ $w\in M,$ and $(w_1-x_3w)w-\frac{1}{2}w=x_2e \in M.$
In both cases, we have  $e\in M$ and $dim(M)=3.$
\end{enumerate}

\subsection{Subalgebras of $D_t(\alpha)$}
The study of subalgebras of the superalgebra $D_t(\alpha)$ has 3 parts. We are considering $1$-dimensional, $2$-dimensional and $3$-dimensional subalgebras.
The result of this subsection is the following

\begin{Th}\label{subd}
Let $M$ be a non-trivial subalgebra of $D_t(\alpha),$ then
\begin{enumerate}

    \item if $dim(M)=1$ and $\alpha\neq\frac{1}{2}, t\neq -1,$ then $M$ is generated by one of the following vectors
     $e_1+e_2,$  $\gamma_1e_i+\gamma_2x,$  $\gamma_1e_i+\gamma_2y,$ where $i=1,2,$ and $\gamma_1,\gamma_2 \in \mathbb{F};$

    \item if $dim(M)=1$ and $\alpha=\frac{1}{2},$ then $M$ is generated by one of the following vectors
     $e_1+e_2,$ $\gamma_1x+ \gamma_2 y,$  $e_1+\gamma_1x+\gamma_2y,$ or  $e_2+\gamma_1y+\gamma_2y,$ where $\gamma_1,\gamma_2 \in \mathbb{F};$

    \item if $dim(M)=1$ and $\alpha\neq\frac{1}{2},t= -1,$ then $M$ is generated by one of the following vectors
    $e_1+e_2,$ $\gamma_1e_i+\gamma_2x,$ $\gamma_1e_i+\gamma_2y,$ where $\gamma_i \in \mathbb{F}, \ i=1,2$ or
    $\gamma_1e_1+\gamma_2e_2+\gamma_3x+\frac{\gamma_1\gamma_2}{\gamma_3(4\alpha-2)}y,$  where  $\gamma_i \in \mathbb{F}^*, i =1,2,3;$


    \item if $dim(M)=2$ and $\alpha\neq\frac{1}{2},t\neq -1,$ then $M$ is generated by one of the following pairs of vectors
    $(e_1,e_2),$ $(e_i, x),$  $(e_i,y),$ $(e_1+e_2, x),$ $(e_1+e_2,y),$ $(e_1+ \gamma x, e_2-\gamma x),$ 
    or $(e_1+ \gamma y, e_2-\gamma y),$ where $\gamma \in \mathbb{F};$ 

    \item if $dim(M)=2$ and $\alpha=\frac{1}{2},t\neq -1,$ then $M$ is generated by one of the following pairs of vectors
    $(e_1,e_2),$ $(e_i, x),$  $(e_i,y),$ $(e_1+e_2, x),$ $(e_1+e_2,y),$ or
    $(e_1+ \gamma_1x+ \gamma_2 y, e_2-\gamma_1 x-\gamma_2 y)$ where $\gamma_1,\gamma_2 \in \mathbb{F};$

    \item if $dim(M)=2$ and $\alpha\neq \frac{1}{2},t=-1,$ then $M$ is generated by one by of the following pairs of vectors
    $(e_1,e_2),$ $(e_i, x),$  $(e_i,y),$   $(e_1+ \gamma_1y, e_2-\gamma_1y),$ $(e_1+ \gamma_1x, e_2-\gamma_1x),$ 
    or $(e_1+e_2, \gamma_1x+ \gamma_2 y),$ where $\gamma_i \in \mathbb{F}, i=1,2;$


    \item if $dim(M)=3$ and $\alpha\neq \frac{1}{2}$ then $M$ is generated one of the following set of vectors
    $(e_1,e_2,x)$ or  $(e_1,e_2,y);$

    \item if $dim(M)=3$ and $\alpha=\frac{1}{2}, t\neq 1$ then $M$ is generated one of the following set of vectors
    $(e_1,e_2,\gamma_1x+\gamma_2y),$ where $\gamma_1,\gamma_2 \in \mathbb{F};$ 

    \item if $dim(M)=3$ and $\alpha=\frac{1}{2}, t= 1$ then $M$ is generated one of the following set of vectors
    $(e_1+e_2, \gamma_1e_2+x,\gamma_2e_2+y)$ or $(e_1,e_2,\gamma_1x+\gamma_2y),$ where $\gamma_1,\gamma_2 \in \mathbb{F}.$

\end{enumerate}

\end{Th}

The proof of the Theorem \ref{subd} is based on the following subsections \ref{1}, \ref{2} and \ref{3}.

\subsubsection{$1$-dimensional subalgebras}\label{1}
Let $M$ be an $1$-dimensional subalgebra of $D_t(\alpha),$ then $M$ is generated by the following vector  
$w_1= x_1e_1+x_2e_2+x_3x+x_4y.$
Now, 
 $$w_1^2=x_1^2e_1+x_2^2e_2+x_1x_3x+x_1x_4y+x_2x_3 x+x_2x_4y+x_3x_4((4\alpha-2)e_1+(2-4\alpha)te_2) \in M.$$
 It is  easy to see that there is some element $k \in \mathbb{F},$ such that 
 $$(x_1+x_2)x_3=kx_3, \ \  x_1^2+x_3x_4(4\alpha -2)=kx_1,$$
 $$(x_1+x_2)x_4=kx_4, \   x_2^2+x_3x_4(2-4\alpha)t=kx_2.$$

Obviously, that if $x_3=x_4=0,$ we have three $1$-dimensional subalgebras. They are subalgebras generated by $e_1, e_2$ and $e_1+e_2.$

If $x_3x_4=0$ and $x_3\neq x_4,$ then $k=x_1+x_2$ and $x_1x_2=0.$
Here we have $4$ opportunities: 
$$x_1=x_3=0, \ x_2=x_3=0, \ x_1=x_4=0 \mbox{ and }x_2=x_4=0.$$
From the first, we have 
$(x_2e_2+x_4y)^2=x_2(x_2e_2+x_4y)$ and we have that $x_2e_2+x_4y$ generates an $1$-dimensional subalgebra.
From other cases, we have that $x_1e_1+x_4y,$ $x_2e_2+x_3x$ and $x_1e_1+x_3x$ are generating $1$-dimensional subalgebras.

On the other hand, if $x_3x_4 \neq 0,$ then $k=x_1+x_2$ and 
$$x_3x_4(4\alpha-2)=x_1x_2= x_3x_4(2-4\alpha)t.$$
Here, if $\alpha=\frac{1}{2}$ we have $x_1x_2=0$ and $M$ is one of subalgebras generated by $e_1+ x_3x+x_4y$ or $e_2+x_3x+x_4y.$
If $\alpha\neq \frac{1}{2},$ it is follows that $t=-1$
and we can see that the superalgebra  $D_{-1}(\alpha)$ has an $1$-dimensional  subalgebra generated by
$x_1e_1+x_2e_2+x_3x+\frac{x_1x_2}{x_3(4\alpha-2)}y,$  where  $x_1, x_2, x_3\neq0.$

\subsubsection{$2$-dimensional subalgebras}\label{2}
Let $M$ be a $2$-dimensional subalgebra of $D_t(\alpha),$ then 
$M$ is generated by the following linear independent vectors
$\alpha_1e_1+\beta_1 e_2+\gamma_1x+\delta_1y$ and $\alpha_2e_1+\beta_2 e_2+\gamma_2x+\delta_2y.$
\begin{enumerate}
    \item[I.] If $\alpha_1e_1+\beta_1 e_2$ and $\alpha_2e_1+\beta_2 e_2$ are linear dependent, we can suppose that 
 $w_1=\alpha_1e_1+\beta_1 e_2+\gamma_1x+\delta_1y$ and $w_2=\gamma_2x+\delta_2y$ are in $M.$
It is easy to see, that $w_2^2=\gamma_2\delta_2(4\alpha-2)(e_1-te_2) \in M.$

\begin{enumerate}

\item[$\bullet$] In the first, let $\gamma_2=0$ and $y\in M,$ then 
$y  w_1+w_1y= (\alpha_1+\beta_1)y + \gamma_1(4\alpha -2)(e_1-te_2) \in M.$

Here, 
if $\gamma_1=0,$ then $M$ is $2$-generated if only  if the element $\alpha_1e_1+\beta_1 e_2$ generates an $1$-dimensional subalgebra.
It is one of the following cases $\alpha_1=0, \beta_1=0$ or $\alpha_1=\beta_1=1.$
In all these three cases, we have examples of $2$-generated subalgebras.
It is one the following subalgebras generated by $(e_1,y), (e_2,y)$ or $(e_1+e_2,y).$ 

If $\gamma_1\neq 0$ and $\alpha\neq\frac{1}{2},$ then $e_1-te_2 \in M$ and $e_1+t^2e_2 \in M.$
If $t\neq -1,$ then $e_1,e_2\in M.$ 
On the other hand, if $t=-1$ we have that $e_1+e_2 \in M$ and $M$ is generated by $e_1+e_2$ and $y.$
    
If $\gamma_1\neq 0$ and $\alpha=\frac{1}{2},$ then 
$(w_1-\delta_1y)^2 = \alpha_1^2e_1+\beta_1^2e_2+(\alpha_1+\beta_1)\gamma_1x \in M$
and for an element $k\in \mathbb{F}$ it is true $\alpha_1^2e_1+\beta_1^2e_2+(\alpha_1+\beta_1)\gamma_1x=k(\alpha_1e_1+\beta_1e_2+\gamma_1x).$
It is follows, that $\alpha_1+\beta_1=k$ and theare are three opportunities: $\alpha_1=k, \beta_1=0,$ or $\alpha_1=0,\beta_1=k,$ or $\alpha_1=\beta_1=0.$
But every this opportunity gives a subalgebra with the dimension $>2.$

\item[$\bullet$] In the second, let $\delta_2=0$ and $x\in M,$ then by some similar way, it is easy to see that
$M$ is generated by $(e_1,x), (e_2,x)$ or $(e_1+e_2,x).$

\item[$\bullet$] In the third, let $e_1-te_2 \in M,$ then $e_1+t^2e_2 \in M.$
If $t\neq-1,$ then $M$ is generated by $e_1$ and $e_2.$
On the other hand, if $t=-1$ we have $e_1+e_2\in M$ and 
$M$ is generated by $e_1+e_2$ and $\gamma_1x+\delta_1y,$ where one or two elements $\gamma_1, \delta_1$ is not equals to $0.$
\end{enumerate}

\item[II.] If $\alpha_1e_1+\beta_1 e_2$ and $\alpha_2e_1+\beta_2 e_2$ are linear independent, then we can suppose that 
 $w_1= e_1+\gamma_1x+\delta_1y$ and $w_2=e_2+\gamma_2x+\delta_2y$ are in $M.$
It is easy to see, that 
$w_i^2=w_i+ \gamma_i\delta_i(4\alpha-2)(e_1-te_2) \in M.$

\begin{enumerate}
    \item[$\bullet$] Suppose that $\alpha=\frac{1}{2}.$
    Here, from $w_1w_2=\frac{1}{2}(\gamma_1+\gamma_2)x+\frac{1}{2}(\delta_1+\delta_2)y \in M,$
    it is easy to see that $\gamma_1+\gamma_2=\delta_1+\delta_2=0.$
    Now we have one family of $2$-dimensional subalgebras of $D_t(\frac{1}{2}),$ generated by 
     $e_1+\gamma_1x+\delta_1y$ and $e_2-\gamma_1x-\delta_1y.$

\item[$\bullet$] Suppose that $\alpha\neq \frac{1}{2}.$ 
\begin{enumerate}
\item[$\star$] In the first, if $\gamma_1=\gamma_2=0,$ then  $w_1  w_2 + w_2  w_1= (\delta_1+\delta_2)y \in M.$

Here, if $\delta_1+\delta_2=0,$ then $M$ is generated by $e_1+ \delta_1y$ and $e_2-\delta_1y.$

If $\delta_1+\delta_2\neq 0,$ then $e_1,e_2, y\in M$ and $M$ is not $2$-generated.

\item[$\star$] In the second, if $\gamma_1=\delta_1=0,$ then  
$w_1 w_2 +w_2 w_1 = \gamma_1\delta_1(4\alpha-2)(e_1-te_2) \in M.$ 
If $w_1 w_2 +w_2 w_1 \neq 0,$ then $w_1 w_2 +w_2 w_1, w_1, w_2$ are linear independent and from here $\gamma_1\delta_1=0.$
It is easy to see, that $M$ is generated by $e_1$ and $e_2.$

\item[$\star$] From the third and 4th cases, by some similar way, we can obtain one new family of $2$-dimensional subalgebras,
generated by $(e_1+ \gamma_1x,e_2-\gamma_1x).$
\end{enumerate}

\end{enumerate}

\end{enumerate}

\subsubsection{$3$-dimensional subalgebras}\label{3}
Let $M$ be a $2$-dimensional subalgebra of $D_t(\alpha),$ then 
$M$ is generated by the following linear independent vectors
$\alpha_ie_1+\beta_i e_2+\gamma_ix+\delta_iy,$ where $i=1,2,3.$
Here we have two opportunities:
\begin{enumerate}
\item[I.] 
If the vector space generated by vectors $\gamma_i x+\delta_i y,$ for $i=1,2,3$ has dimension $2,$ then 
we can suppose that 
  $w_1=\alpha_1e_1+\beta_1 e_2+x,$
  $w_2=\alpha_2e_1+\beta_2 e_2+y$ and $w_3=\alpha_3e_1+\beta_3 e_2$ are in $M.$
It is easy to see, that $w_3^2=\alpha_3^2 e_1+ \beta_3^2e_2 \in M.$
If $w_3^2$ and $w_3$ are linearly independent, then $e_1,e_2 \in M$ and $dim(M)=4.$
On the other hand, 
we can consider $w_3$ as $e_1, e_2$ or $e_1+e_2.$ 
\begin{enumerate}

\item[$\bullet$] From the first, it is easy to see that we can consider $w_1$ and $w_2$ as 
$\beta_1 e_2+x$ and $\beta_2 e_2+y.$
Now, $e_1(\beta_1 e_2+x) , (\beta_2 e_2+y)e_1 \in M$ and $e_1,x,y, xy \in M, dim(M)=4.$

\item[$\bullet$] From the second, we have a similar result.

\item[$\bullet$] From the third opportunity, we have $e_1+e_2 \in M$
and we can consider $w_1$ and $w_2$ as 
$\beta_1 e_2+x$ and $\beta_2 e_2+y.$
Here, $(\beta_2x-\beta_1y)^2= \beta_1\beta_2(4\alpha-2)(e_1-te_2) \in M.$
\begin{enumerate}
\item[$\star$] $\alpha=\frac{1}{2},$ then $(\beta_1e_2+x)(\beta_2e_2+y) \in M$ and $t=1.$
The subalgebra generated by $e_1+e_2, \beta_1e_2+x, \beta_2e_2+y$ is a $3$-dimensional subalgebra of $D_1(\alpha).$

\item[$\star$] If $\beta_1=0$ or $\beta_2=0,$ then $xy, yx \in M$ and $\alpha=(1-\alpha)t, 1-\alpha=\alpha t.$ 
It is follow that $\alpha = \alpha t^2$ and $t=\pm 1,$ $\alpha=\frac{1}{2}$ and we have the first case.

\item[$\star$] $t=-1,$ then $(\beta_1e_2+x)(\beta_2e_2+y), (\beta_2e_2+y)(\beta_1e_2+x) \in M$ and
$\beta_1\beta_2(1-2\alpha)+2=0, \beta_1\beta_2(1-2(1-\alpha))+2=0.$
There are no elements $\beta_1, \beta_2, \alpha$ which satisfying this condition.
\end{enumerate}
\end{enumerate}

\item[II.] 
If the vector space generated by vectors $\gamma_i x+\delta_i y,$ for $i=1,2,3$ has dimension $1,$ then 
we can suppose that 
 $w_1=\alpha_1e_1+\beta_1 e_2+\gamma_1 x+\delta_1 y,$
  $w_2=\alpha_2e_1+\beta_2 e_2$ and $w_3=\alpha_3e_1+\beta_3 e_2$ are in $M.$
Obviously, $w_2$ and $w_3$ are linearly independent and $e_1,e_2 \in M.$
Now, $\gamma_1 x+\delta_1 y \in M$ and 
$e_1(\gamma_1 x+\delta_1 y), (\gamma_1 x+\delta_1 y)e_1$ are linearly dependent. 
Now, for $M$ we have only two opportunities:
\begin{enumerate}
    
\item[$\bullet$] $\alpha\neq \frac{1}{2}$ and $M$ is generated by $(e_1,e_2,x)$ or $(e_1,e_2,y).$

\item[$\bullet$] $\alpha= \frac{1}{2}$ and $M$ is generated by $(e_1,e_2, \gamma_1x+\delta_1y).$

\end{enumerate}

\end{enumerate}

\subsection{Subalgebras of $D_t^{1/2}$}
The study of subalgebras of the superalgebra $D_t^{1/2}$ has 3 parts. 
We are considering $1$-dimensional, $2$-dimensional and $3$-dimensional subalgebras.
The result of this subsection is the following

\begin{Th}\label{subd12}
Let $M$ be a non-trivial subalgebra of $D_t^{1/2},$ then
\begin{enumerate}
\item if $dim(M)=1$ and $t\neq -1,$ then $M$ is generated by one of the following vectors $e_1+e_2,$ $\gamma_1e_1+\gamma_2y$ or $\gamma_1e_2+\gamma_2y,$ where $\gamma_i \in \mathbb{F};$

\item if $dim(M)=1$ and $t= -1,$ then $M$ is generated by $e_1+e_2,$ or $\gamma_1e_1+\gamma_2e_2+\sqrt{-\gamma_1\gamma_2}x+\gamma_3y,$ where $\gamma_i \in \mathbb{F};$


\item if $dim(M)=2$ and $t\neq-1,$ then 
$M$ is generated by one of the following pairs of vectors $(e_i, y), i=1,2,$ or $(e_1+\gamma y, e_2-\gamma y),$ where $\gamma \in \mathbb{F};$

\item if $dim(M)=2$ and $t=-1,$ then 
$M$ is generated by one of the following pairs of vectors 
$(e_i, y), i=1,2,$ $(e_1+\gamma y, e_2-\gamma y),$  or $(e_1+e_2, \gamma_1 x+\gamma_2 y),$ 
where $\gamma_1, \gamma_2 \in \mathbb{F};$


    \item if $dim(M)=3$ and $t\neq-1,$ then  $M$ is generated of  $(e_1,e_2,y);$

    \item if $dim(M)=3$ and $t=-1,$ then  $M$ is generated one of the following sets of vectors 
    $(e_1,e_2,y)$ or $(e_1+e_2, \gamma e_2+x, \pm 2i e_2 +y),$ where $\gamma \in \mathbb{F}.$

\end{enumerate}
\end{Th}

The proof of the Theorem \ref{subd12} is based on the following subsections \ref{dt1}, \ref{dt2} and \ref{dt3}.

\subsubsection{$1$-dimensional subalgebras}\label{dt1}
Let $M$ be an $1$-dimensional subalgebra, 
then $M$ is generated by the following vector $w_1=\alpha_1e_1+\alpha_2e_2+\alpha_3 x+\alpha_4 y$
and there is some element $k\in \mathbb{F},$ such that
$$kw_1=w_1^2=\alpha_1^2e_1+\alpha_2^2e_2+\alpha_1\alpha_3x+\alpha_1\alpha_4y+\alpha_2\alpha_3x+\alpha_2\alpha_4y-\alpha_3^2(e_1-te_2).$$
Here, 
$$k\alpha_1=\alpha_1^2-\alpha_3^2, \ k\alpha_2=\alpha_2^2+t\alpha_3^2, \ k\alpha_3=(\alpha_1+\alpha_2)\alpha_3, \ k\alpha_4=(\alpha_1+\alpha_2)\alpha_4.$$

If $x_3=x_4=0,$ there are three $1$-dimensional subalgebras generated by $e_1, e_2,$ or $e_1+e_2.$

If $x_3\neq0$ or $x_4\neq0,$ then $k=\alpha_1+\alpha_2$ and $t=-1$ or $\alpha_3=\alpha_1\alpha_2=0.$

\begin{enumerate}
\item[$\bullet$] In the first opportunity, we have that $M$ is generated by $\alpha_1e_1+\alpha_2e_2+\sqrt{-\alpha_1\alpha_2}x+\alpha_4y.$

\item[$\bullet$] In the second opportunity, we have that $M$ is generated by one of the following vectors $\alpha_1e_1+\alpha_4y$ or $\alpha_2e_2+\alpha_4y.$
\end{enumerate}

\subsubsection{$2$-dimensional subalgebras}\label{dt2}
 Let $M$ be a $2$-dimensional subalgebra, 
then $M$ is generated by the following linear independent vectors 
$w_i=\alpha_ie_1+ \beta_ie_2+\gamma_i x+\delta_i y, i=1,2.$
Then

\begin{enumerate}
    
    \item[I.] If $\alpha_ie_1+ \beta_ie_2, i=1,2$ are linear dependent, we can consider $w_1$ and $w_2$ as 
    $\alpha_1e_1+ \beta_1e_2+\gamma_1 x+\delta_1y$ and $\gamma_2 x+\delta_2 y.$
    Now, it is easy to see, that $w_2^2=-\gamma_2^2(e_1-te_2) \in M.$ 
    Here, if $\gamma_2 \neq 0,$ then $t=-1$ and $M$ is  generated by $e_1+e_2$ and $\gamma_2 x+\delta_2 y.$
    On the other hand,  $\gamma_2=0,$ and  we have $y\in M$ and $w_1y-yw_1=2\gamma_1(e_1+te_2) \in M.$
    Obvioulsy, 
    
\begin{enumerate}

\item[$\bullet$] if $\gamma_1=0,$ then $\alpha_1 e_1+\alpha_2e_2, y \in M$ and $M$ is generated by one of the following pairs of vectors $(e_1+e_2,y), (e_i, y), i=1,2.$
\item[$\bullet$] If $\gamma_1\neq0,$ then $e_1+te_2,x,y\in M$ and $dim(M)>2.$
\end{enumerate}        
    
    \item[II.] If $\alpha_ie_1+ \beta_ie_2, i=1,2$ are linear independent, we can consider $w_1$ and $w_2$ as 
    $e_1\\+\gamma_1 x+\delta_1y$ and $e_2+\gamma_2 x+\delta_2 y.$
    Now, it is easy to see, that $w_i^2=w_i-\gamma_i^2(e_1-te_2), i=1,2.$ 
    
    \begin{enumerate}
    
    \item[$\bullet$] If $\gamma_i=0 \ (i=1,2),$ then      $w_1w_2=\frac{1}{2}(\delta_1+\delta_2)y \in M$ and $\delta_1+\delta_2=0.$
    Here, $M$ is generated by $e_1+\delta_1y$ and $e_2-\delta_1 y.$
    
    \item[$\bullet$] If $\gamma_1\neq0$ or $\gamma_2\neq0,$ then $e_1-te_2\in M$ and $e_1+t^2e_2\in M.$
    It is follows that $t=-1, e_1+e_2\in M$ or $dim(M)>2.$
    If $t=-1,$ then $M$ is generated by $e_1+e_2$ and $\gamma_1x+\delta_1y.$
    \end{enumerate}

\end{enumerate}

\subsubsection{$3$-dimensional subalgebras}\label{dt3}
Let $M$ be a $3$-dimensional subalgebra, 
then $M$ is generated by the following linear independent vectors 
$w_i=\alpha_ie_1+ \beta_ie_2+\gamma_i x+\delta_i y, i=1,2,3.$
Then

\begin{enumerate}
    
    \item[I.] If $\gamma_ix+ \delta_iy, i=1,2,3$ are generating a subspace of dimension $2,$ then 
    we can consider $w_1=\alpha_1e_1+ \beta_1e_2+x, w_2=\alpha_2e_1+ \beta_2e_2+ y, w_3=\alpha_3e_1+ \beta_3e_2.$ 
    \begin{enumerate} 
    \item[$\bullet$] If $w_3^2$ and $w_3$ are linear independent, then $dim(M)=4.$
    \item[$\bullet$] If $w_3^2$ and $w_3$ are linear dependent, then $w_3$ can be considered as $e_1,e_2$ or $e_1+e_2.$
    \begin{enumerate}
    \item[$\star$] If $w_3=e_1,$ then $e_1(\beta_1e_2+x)$ and $(\beta_2e_2+x)e_1$ are linearly independent and $dim(M)=4.$
    \item[$\star$] If $w_3=e_2,$ then also $dim(M)=4.$
    \item[$\star$]  If $w_3=e_1+e_2,$ then 
    we can consider $w_1=\beta_1e_2+x$ and $w_2=\beta_2e_2+y,$ then 
    $w_1^2-\beta_1w_1$ and $e_1+e_2$ are linearly dependent and $t=-1.$
    Also, note that 
    $2w_1w_2=\beta_1 w_2+\beta_2w_1+2e_1+(-2-\beta_2^2)e_2 \in M$ and $\beta_2 =\pm 2i.$
    It is easy to see, that $M$ is generated by $e_1+e_2, \beta_1e_2+x$ and $ \pm 2i e_2 +y.$
    \end{enumerate}
    \end{enumerate}

    \item[II.] If $\gamma_ix+ \delta_iy, i=1,2,3$ are generating a subspace of dimension $1,$ then 
    we can consider $w_1=\gamma_1x+\delta_1 y, w_2=e_1, w_3=e_2.$
    It is easy to see that $e_1w_1$ and $e_2w_2$ are linearly dependent and $\gamma_1=0$
    and $M$ is generated by $e_1,e_2,$ and $y.$

    \end{enumerate}


\section{Automorphisms of $K_3(\alpha,\beta,\gamma)$ and $D_t(\alpha,\beta,\gamma)$}

The main idea of the description of automorhisms is following:
if $e^2=e$ and $x^2=0,$ then for any automorphism $\phi$ we have
$\phi(e)^2=\phi(e)$ and $\phi(x)^2=0.$
From the classification of $1$-dimensional subalgebras of $K_3(\alpha,\beta,\gamma)$ and $D_t(\alpha,\beta,\gamma)$
we can obtain images of elements with this property.
Also note that every automorphism of the algebra $A$ is an automorphism of the algebra $A^{(+)}.$

\subsection{Automorphisms of $K_3(\alpha)$}
Let $\phi$ be an automorphism of $K_3(\alpha)$ 
and  $\phi(x)=\phi_{x,z}z+\phi_{x,w}w,$ where $x \in \{z,w\}.$
Then,  from the Theorem \ref{subk} we have the following two opportunities: 

\begin{enumerate}
    \item[I.] $\alpha\neq \frac{1}{2}$ and  $\phi(e)=e+kz+lw, \phi(z)=\phi_{z,z}z, \phi(w)=\phi_{w,w}w.$ Then 
    $$2\phi(e)=\phi(zw-wz)=2\phi_{z,z}\phi_{w,w}e.$$ 
    In the end, it is easy to verify that the linear mapping $\phi$ defined as
    \begin{eqnarray}\label{autk3_1}\phi(e)=e, \ \phi(z)= \gamma z, \ \phi(w) = \gamma^{-1} w, \mbox{ where } \gamma \in \mathbb{F}
    \end{eqnarray}
    is an automorphism.

    \item[II.] $\alpha=\frac{1}{2}$ and $\phi(e)=e+kz+lw.$ 
    Then,
        $$e+kz+lw=\phi(e)=\phi(z)\phi(w)=(\phi_{z,z}\phi_{w,w}-\phi_{z,w}\phi_{w,z})e$$
    and $\phi(e)=e.$
    In the end, it is easy to verify that the linear mapping $\phi$ defined as
    \begin{eqnarray}\label{autk3}\phi(e)=e, \ \phi(z)= \gamma_1 z+\gamma_2 w, \ \phi(w) = \gamma_3z+\gamma_4 w, \mbox{ where } \gamma_i \in \mathbb{F} \mbox{ and } \gamma_1\gamma_4-\gamma_2\gamma_3=1.
    \end{eqnarray}
    is an automorphism.

\end{enumerate}

\subsection{Automorphisms of $K_3^{1/2}$}

Let $\phi$ be an automorphism of  $K_3^{1/2}$.
Then it is easy to see, that  
$\phi(e)=e+x_3w$ and $\phi(z)=\phi_{z,z}z+\phi_{z,w}w,\phi(w)=\phi_{w,w}w.$
Now,
it is easy to see, that $\phi(e)=-\phi(z)\phi(z)=\phi_{z,z}^2e$ and $\phi(e)=e.$
Also, from $\phi(e)=\phi(zw)=\phi_{z,z}\phi_{w,w}e$ we obtain that $\phi_{z,z}=\phi_{w,w}.$

In the end, we can verify that the linear mapping
    \begin{eqnarray} \label{autk12}\phi(e)=e, \ \phi(z)=\pm z+kw, \ \phi(w)=\pm w.    \end{eqnarray}
gives an automorphism for every element $k \in \mathbb{F}.$

\subsection{Automorphisms of $D_{t}(\alpha)$}
Let $\phi$ be an automorphism of $D_{t}(\alpha)$. 
Then it is easy to see, that 
$$\phi(e_1)=e_1+\gamma x+ \delta y  \mbox{ and }\phi(e_2)=e_2-\gamma x - \delta y.$$
Then,  from the Theorem \ref{subd} we have the following two opportunities: 

\begin{enumerate}

    \item[I.] $\alpha\neq \frac{1}{2}$ and $\phi(x)=\phi_{x,x}x$ and $\phi(y)=\phi_{y,y}y.$
Now, 
$$2\alpha\phi(e_1)+(1-\alpha)t\phi(e_2)=\phi(xy)=\phi_{x,x}\phi_{y,y}(\alpha e_1+(1-\alpha)te_2)$$
and $\phi_{x,x}=\phi_{y,y}^{-1}$ and $\gamma=\delta=0.$
In the end, it is easy to verify that the linear mapping $\phi$ defined as
\begin{eqnarray}
\phi(e_1)=e_1, \ \phi(e_2)=e_2, \ \phi(x)= \gamma x, \ \phi(y)=\gamma^{-1}y
\end{eqnarray}
is an automorphism.

    \item[II.] $\alpha=\frac{1}{2}$ and $\phi(x)=\gamma_1 x+\gamma_2y, \phi(y)=\gamma_3 x+\gamma_4y.$
Now, 
$$\phi(e_1+te_2)=\phi(xy)=(\gamma_1\gamma_4-\gamma_2\gamma_3)(e_1+te_2)$$
and $\gamma_1\gamma_4-\gamma_2\gamma_3=1$ and $\gamma=\delta=0.$
In the end, it is easy to verify that the linear mapping $\phi$ defined as
\begin{eqnarray}
\phi(e_1)=e_1, \ \phi(e_2)=e_2, \ \phi(x)= \gamma_1 x+\gamma_2y, \ \phi(y)=\gamma_3x+\gamma_4y, 
\mbox{ where } \gamma_1\gamma_4-\gamma_2\gamma_3=1\end{eqnarray}
is an automorphism.

\end{enumerate}

\subsection{Automorphisms of $D_{t}^{1/2}$}
Let $\phi$ be an automorphism of $D_{t}^{1/2}$.
Then it is easy to see, that  
$$\phi(e_i)=e_i,\phi(x)=\phi_{x,x}x+\phi_{x,y}y, \phi(y)= \phi_{y,y}y.$$

    It is easy to see, that 
    $\phi(x+y)=2\phi(e_1x)=\phi_{x,x}x+\phi_{x,x}y+\phi_{x,y}y$
    and $\phi_{x,x}=\phi_{y,y}.$
    Obviously, 
    $$e_1+te_2=\phi(e_1+te_2)=\phi(x)\phi(y)=\phi_{x,x}^2(e_1+te_2)$$
    and $\phi_{x,x}=\phi_{y,y}=\pm1.$
    In the end, we can verify that every linear mapping defined as
\begin{eqnarray}
\phi(e_1)=e_1, \ \phi(e_2)=e_2, \ \phi(x)= \pm x+ \gamma y, \ \phi(y)= \pm y
\end{eqnarray}
is an automorphism.

\section{Derivations of the simple noncommutative Jordan superalgebras}

\medskip
Our basic principle of determining the derivations of a simple noncommutative Jordan superalgebra $U$ is the following. Firstly, we compute the superalgebra of derivations of the symmetrized superalgebra $U^{(+)},$ which structure was specified earlier. Then we check which derivations of $U^{(+)}$ are derivations of $U.$ We will do this by verifying (\ref{der_con}) on the product of basis elements of $U^{(+)}$ and $U$. Further on, when we write down system of relations arising from (\ref{der_con}), it will be assumed that we write down all nontrivial relations, that is, a mapping $d \in End(U)$ is a derivation of $U$ if and only if it satisfies these conditions. If we know the exact expression for a Poisson bracket $[,]$ on $U^{(+)}$ such that $ab = a\cdot b + [a,b]$ then $Der(U) = Der(U^{(+)}) \cap Der(U,[,]).$

\subsection{Derivations of $K_3(\alpha, \beta, \gamma)$}
By some trivial computations we can obtain the following

\begin{Th}
Let $\mathbb{K}$ be a superalgebra of type $K_3(\alpha,\beta,\gamma).$ 
Then

\begin{enumerate}
    
\item if $\mathbb{K}=K_3,$ then 

$Der(\mathbb{K})_{\bar{0}} \cong sl_2,$  and even derivation $d$ act by 
$$ed=0, \ zd= \gamma_1z + \gamma_2 w, \ wd=  \gamma_3 z - \gamma_1 w, \mbox{ where }\gamma_1,\gamma_2,\gamma_3 \in \mathbb{F};$$

$Der(\mathbb{K})_{\bar{1}}$ is $2$-dimensional, and odd derivation $d$ act by 
$$ed=\gamma_1 z + \gamma_2 w, \ zd= -2\gamma_2e, \ wd = 2\gamma_1e,  \mbox{ where }\gamma_1,\gamma_2 \in \mathbb{F}; $$

\item if $\mathbb{K}=K_3(\alpha), \alpha \neq \frac{1}{2},$ then 
$Der(\mathbb{K})$ is a one-dimensional Lie algebra,

and even derivation $d$ act by 
$$ed=0, \ zd= \gamma z, \ wd=   - \gamma w, \mbox{ where }\gamma \in \mathbb{F};$$

\item if $\mathbb{K}=K_3^{1/2},$ then 
$Der(\mathbb{K})$ is a one-dimensional Lie algebra,

and even derivation $d$ act by 
$$ed=0, \ zd= 0, \ wd=   \gamma w, \mbox{ where }\gamma \in \mathbb{F}.$$
    
\end{enumerate}

\end{Th}

\subsection{Derivations of $D_t(\alpha)$}
By some trivial computations we can obtain the following

\begin{Th}
Let $\mathbb{K}$ be a superalgebra of type $D_t(\alpha,\beta,\gamma).$ 
Then

\begin{enumerate}
    
\item if $\mathbb{K}=D_t,$ then 

$Der(\mathbb{K})_{\bar{0}} \cong sl_2,$  and even derivation $d$ act by 
$$e_1d=0, \ e_2d=0, \ xd= \gamma_1x + \gamma_2 y, \ yd=  \gamma_3 x - \gamma_1 y, \mbox{ where }\gamma_1,\gamma_2,\gamma_3 \in \mathbb{F};$$

$Der(\mathbb{K})_{\bar{1}}$ is $2$-dimensional, and odd derivation $d$ act by 
$$e_1d=\gamma_1 x + \gamma_2 y, \ e_2d=-(\gamma_1x+\gamma_2 y), \ xd= 2(-\gamma_2 e_1+ \gamma_2te_2), \ wd = 2(\gamma_1e_1-\gamma_1te_2),  \mbox{ where }\gamma_1,\gamma_2 \in \mathbb{F}; $$

\item if $\mathbb{K}=D_t(\alpha), \alpha \neq \frac{1}{2},$ then 
$Der(\mathbb{K})$ is a one-dimensional Lie algebra,

and even derivation $d$ act by 
$$e_1d=0, \ e_2d=0, \ xd= \gamma x, \ yd=   - \gamma y, \mbox{ where }\gamma \in \mathbb{F};$$

\item if $\mathbb{K}=D_t^{1/2},$ then 
$Der(\mathbb{K})$ is a one-dimensional Lie algebra,

and even derivation $d$ act by 
$$e_1d=0, \ e_2d=0, \ xd= 0, \ yd=   \gamma y, \mbox{ where }\gamma \in \mathbb{F}.$$
    
\end{enumerate}

\end{Th}

\subsection{Derivations of $U(V, f, \star)$ }
At first we calculate the derivations of $J(V,f).$

\subsubsection{Even derivations of $J(V,f)$}
Let $d$ be an even derivation of $J(V,f), 0 \neq v \in V_0$, and let $vd = \alpha + u,$ where $\alpha \in \mathbb{F}, u \in V_0$. Then $0 = v^2d = 2v\cdot vd = 2 \alpha v + 2 f(u,v)$, therefore $\alpha = 0,$ $f(v,vd) = 0$. It is easy to see that this relation and its linearization $f(wd,v) - f(w,vd) = 0$ hold also for $w, v \in V_1$, and also that they are sufficient for $d \in End(J(V,f))$ to be an even derivation of $J(V,f).$ Hence,
\begin{equation}
\label{JVf_even}
Der(J(V,f))_{\bar{0}} = \{d \in End(V)_{\bar{0}} : f(wd,v) - f(w,vd) = 0, \text{where } w, v \in V_0 \cup V_1\}.
\end{equation}
\subsubsection{Odd derivations of $J(V,f)$}
Suppose that $V_0 \neq 0$. Let $v \in V_0, w \in V_1,$ and $vd = u \in V_1, wd = \alpha + t, \alpha \in \mathbb{F}, t \in V_0.$ Then $0 = wvd = (\alpha + t)v + f(w,u)$, and again we have $\alpha = 0$ and $V_1d \subseteq V_0.$ We infer that $f(wd,v) = -f(w,vd)$ for all $w \in V_1, v \in V_0.$ It is easy to see that this relation is sufficient for $d$ to be an odd derivation of $J(V,f)$. Hence,
\begin{equation}
\label{JVf_odd}
Der(J(V,f))_{\bar{1}} = \{d \in End(V)_{\bar{1}} : f(wd,v) + f(w,vd) = 0, \text{where } w \in V_1, v \in V_0\}.
\end{equation}
Combining (\ref{JVf_even}) and (\ref{JVf_odd}), we see that, for $s = \bar{0}, \bar{1},$
$$Der(J(V,f))_s = \{d \in End(V)_s : f(wd,v) = (-1)^{s\cdot w}f(w,vd), \text{where } w, v \in V_0 \cup V_1\}.$$
The vector space $Der(J(V,f))$ with respect to the commutator forms a Lie superalgebra which is denoted by $Lieosp(V, f)$.
Now, it is easy to see that the linear mapping $d$ is a derivation of $U(V, f, \star)$ is a derivation of $U$ if and only if it is a derivation of $U^{(+)} = J(V,f)$ and a derivation of superanticommutative product $\star$. So, we have the following result:
\begin{proposition}
$Der(U(V, f, \star)) = Lieosp(V, f) \cap Der(V, \star).$
\end{proposition}
\medskip

\subsection{Derivations of $J(\Gamma_n, A)$}
In this section we consider $\Gamma_n$ as a \textit{Poisson superalgebra}, that is, a superalgebra with two multiplications $\cdot$ and $\{,\}$, where $\{,\}$ is a Poisson bracket defined by (\ref{Pois_br}). By derivation of $\Gamma_n(\cdot, \{,\})$ we mean the derivation with respect to both products. Since $(\Gamma_n, \cdot)$ is a free associative-supercommutative superalgebra with odd generators $x_1, \dots, x_n$, every derivation of $(\Gamma_n, \cdot)$ is uniquely determined by its actions on $x_1, \dots, x_n.$ Particularly, let $d \in Der((\Gamma_n, \cdot))$, and let $x_1d = f_1, \dots, x_nd = f_n,$ where $f_i \in \Gamma_n.$ Then $d = \sum_{i=1}^n \partial_i f_i$. The superalgebra of all derivations of $(\Gamma_n, \cdot)$ is denoted by $W_n.$ The derivations with respect to the both products form the simple Lie subalgebra $H_n$ of $W_n$ \cite{Ret} defined by the following relations:
$$f_i \partial_j = f_j \partial_i, \text{where } i, j = 1, \dots, n.$$
By \cite{Kac}, for every $d = \sum_{i = 1}^n f_i \partial_i \in H_n$ there exists $f \in \Gamma_n$ such that $f_i = f \partial_i,$ so every derivation of this superalgebra has the form $\sum_{i=1}^n(f \partial_i) \partial_i, f \in \Gamma_n.$

The derivation superalgebra of $J(\Gamma_n)$ was computed by Retakh in \cite{Ret}. It was proven that every derivation of $J(\Gamma_n)$ is the sum of the derivations of the two following types:

\begin{enumerate}
\item[D$1)$] Let $d$ be a derivation of $(\Gamma_n, \cdot, \{,\})$. Extend its action on $J(\Gamma_n)$ by setting $\bar{a}d = (-1)^d \overline{ad}.$ It is easy to check that the resulting mapping is indeed a derivation of $J(\Gamma_n)$ and is uniquely defined by the initial derivation of $\Gamma_n.$

\item[D$2)$] Let $x$ be an element of $\Gamma_n.$ Define the linear mapping $0^x$ by $a0^x = 0, \bar{a}0^x = ax,$ where $a \in \Gamma_n.$ One can check directly that $0^x$ is a derivation of $J(\Gamma_n).$
\end{enumerate}

By Retakh's results every derivation of $J(\Gamma_n)$ is uniquely defined by a derivation $d \in Der(\Gamma_n)$ and an element $x \in \Gamma_n.$ We will denote it by $d^x.$ So, the superalgebra of derivations of $J(\Gamma_n)$ is a semidirect sum of $H(n)$ and $\Gamma_n.$
We need to check which derivations of $J(\Gamma_n)$ are derivations with respect to the Poisson bracket in this superalgebra.

\begin{enumerate}

\item Let $a, b \in (\Gamma_n)_{\bar{0}} \cup (\Gamma_n)_{\bar{1}}.$ Then
$$\{\bar{a},\bar{b}\}d = (-1)^{(b+1)d} \{\bar{a}d, \bar{b}\} + \{\bar{a},\bar{b}d\} = (-1)^{bd}\{\bar{ad},\bar{b}\} + (-1)^d\{\bar{a},\bar{bd}\} = $$
$$= (-1)^{bd + b}ad\cdot bA + (-1)^b a\cdot bd A = (-1)^b(abd \cdot A)$$
On the other hand, we have
$$\{\bar{a},\bar{b}\}d = (-1)^b (ab\cdot A)d = (-1)^b (abd \cdot A + ab\cdot Ad),$$
hence $d$ is a derivation of $J(\Gamma_n, A)$ of type D$1)$ if and only if $ab \cdot Ad = 0$ for all homogeneous $a, b \in \Gamma_n$, which is obviously equivalent to $Ad = 0.$

\item  Let $x$ be a homogeneous element of $\Gamma_n$ and let $0^x$ be a derivation of $J(\Gamma_n, A)$. Then
$$0 = (-1)^b abA0^x = \{\bar{a}, \bar{b}\}0^x = (-1)^{(b+1)(x+1)}\{\bar{a}0^x, \bar{b}\} + \{\bar{a}, \bar{b}0^x\} = $$
$$ = (-1)^{(b+1)(x+1)}\{ax,\bar{b}\} + \{\bar{a},bx\} = 0,$$
so every derivation of type D$2)$ is a derivation of $J(\Gamma_n, A)$.
\end{enumerate}

Combining the above results, we have
\begin{Th}
$Der(J(\Gamma_n)) = (H(n)\cap Ann(A)) \rtimes \Gamma_n.$
\end{Th}

\subsection{Derivations of $\Gamma_n(\mathcal{D})$}
From the previous section we know that every derivation of $\Gamma_n$ is of the form
\begin{equation}
\label{der_gr}
d = \sum_{k=1}^n \partial_k f_k, f_k \in \Gamma_n.
\end{equation}
 We need to check which homogeneous derivations of $(\Gamma_n,\cdot)$ are derivations with respect to the Poisson bracket $[,]$.
 Let $d$ be a homogeneous derivation of $(\Gamma_n, \cdot)$ given by (\ref{der_gr}). Then for any homogeneous $h \in \Gamma_n$ we have
 $$ (h\partial_i d - (-1)^d (hd) \partial_i) =  \sum_k(h\partial_i \partial_k \cdot f_k - (-1)^d(h \partial_k \cdot f_k) \partial_i = $$
 $$ = \sum_k(h\partial_i \partial_k \cdot f_k + h\partial_k \partial_i \cdot f_k - (-1)^d h\partial_i \cdot f_k \partial_i)= $$
 $$ = -(-1)^d\sum_k h\partial_k \cdot f_k\partial_i.$$
 Hence for arbitrary two homogeneous $f, g \in \Gamma_n$ we have
 $$ (fg)d - (-1)^{gd}(fd)g - f(gd) = (f\cdot g)d - (-1)^{gd}(fd)\cdot g - f\cdot(gd) +$$
 $$+(-1)^{g+1} \sum_{i,j} (f \partial_i \cdot g \partial_j \cdot a_{ij})d - (-1)^{gd} (fd)\partial_i \cdot g \partial_j \cdot a_{ij} -(-1)^d f \partial_i \cdot (gd)\partial_j \cdot a_{ij}) =$$
 $$=(-1)^{g+1} \sum_{i,j} (-1)^{(g+1)d}(f \partial_i d -(fd)\partial_i) \cdot g \partial_j \cdot a_{ij} + f \partial_i \cdot (g \partial_j d - (gd)\partial_j)\cdot a_{ij} + f \partial_i \cdot g \partial_j \cdot a_{ij}d =$$
 $$= (-1)^{g+1}\sum_{i,j,k}-(-1)^{gd} f\partial_k  \cdot f_k \partial_i \cdot g \partial_j \cdot a_{ij} -(-1)^d f \partial_i \cdot g \partial_k \cdot f_k \partial_j \cdot a_{ij} + f \partial_i \cdot g \partial_j \cdot a_{ij} \partial_k \cdot f_k=$$
 $$ = (-1)^{g+1} \sum_{i,j,k} -(-1)^d f \partial_i \cdot g \partial_j \cdot f_i \partial_j \cdot a_{jk} -(-1)^d f \partial_i \cdot g \partial_j \cdot f_j \partial_k \cdot a_{ik} + f \partial_i \cdot g \partial_j \cdot a_{ij} \partial_k \cdot f_k=$$
 $$ =(-1)^{g+1} \sum_{i,j} f \partial_i \cdot g \partial_j \cdot \Big( \sum_k -(-1)^d f_i \partial_j \cdot a_{jk} - (-1)^d f_j \partial_k \cdot a_{ik} +  a_{ij} \partial_k \cdot f_k\Big).$$

Hence, we have proved the following 
\begin{proposition}
\emph{A homogeneous mapping $d$  of the form (\ref{der_gr}) is in $ Der(\Gamma_n(\mathcal{D}))$ if and only if}
\begin{equation}
\label{gras_der}
\sum_{k=1}^n a_{ij} \partial_k f_k = (-1)^d(\sum_{k=1}^nf_i\partial_k\cdot a_{jk} + f_j\partial_k \cdot a_{ik}), i, j \leq n.
\end{equation}
\end{proposition}
We provide some examples:

$1)$ Let $d_i = \partial_i,$ so the matrix $(a_{ij})_{i, j = 1}^n$ is an identity matrix of size $n.$
Then the (\ref{gras_der}) takes the form
$$0 = f_i \partial_j + f_j \partial_i, i, j \leq n.$$
As in the previous section, such mappings form a simple Lie subalgebra $H_n$ in $W_n.$ 

$2)$ Let $n = 3$ and
$$(a_{ij}) = A = \begin{pmatrix} 1 & 0 & 0\\ 0 & 1 & 0 \\ 0 & 0 & x_1x_2 \end{pmatrix}$$
Then (\ref{gras_der}) takes the following form:
$$ \begin{cases}
f_1 \partial_1 = 0,\\
f_1 \partial_2 + f_2 \partial_1 = 0,\\
f_3 \partial_1 + f_1 \partial_3 x_1 x_2 = 0,\\
f_2 \partial_2 = 0,\\
f_3 \partial_2 + f_2 \partial_3 x_1 x_2 = 0,\\
x_1 f_1 - x_2 f_2 = (-1)^d f_3 \partial_3.
\end{cases} $$
Consider the second equation in this system. Combining the facts that $f_3 \partial_1$ is independent on $x_1$, and $f_1 \partial_3 x_1 x_2$ depends on $x_1$, and both elements lack constant term, the only possibility is that $f_3 \partial_1 = f_1 \partial_3 x_1 x_2 = 0.$ Analogously, $f_3 \partial_2 = f_2 \partial_3 x_1 x_2 = 0.$ Therefore, $f_3$ depends only on $x_3$, and the right part in the last equation is a constant, but the left part has no constant terms, so both sides are equal to zero, therefore $f_3 = \gamma \in \mathbb{F},$ and $x_1f_1 = x_2f_2.$ This implies that both $f_1, f_2$ have zero constant terms. Since $f_1 \partial_1 = 0 = f_2 \partial_2; f_1 = \alpha x_2$ and $f_2 = -\alpha x_1$ if $d$ is even; $f_1 = \beta x_2x_3$ and $f_2 = -\beta x_1x_3$ if $d$ is odd. 
Hence an arbitrary even derivation $d$ of this superalgebra is defined by the following way 
$$x_1d = \gamma_1 x_2, \ x_2d= -\gamma_1 x_1, \ x_3d= \gamma_2, \ \mbox{where} \ \gamma_1, \gamma_2 \in \mathbb{F},$$ 
and an arbitrary odd derivation $d$ of this superalgebra is defined by the following way
$$x_1d=  \gamma  x_2x_3, \ x_2d= -\gamma  x_1x_3, \ x_3d= 0, \ \mbox{where} \ \gamma \in \mathbb{F}.$$

We can also give some partial solutions of (\ref{gras_der}) for arbitrary system of derivations $\mathcal{D}$. We will rewrite (\ref{gras_der}) in a different form. Again, we start with
$$a_{ij}d = [x_i,x_j]d = (-1)^d[x_id, x_j] + [x_i,x_jd] = (-1)^d(x_idd_j + x_jdd_i),$$
so
$$x_id_jd - (-1)^d(x_idd_j + x_jdd_i) = 0.$$
Adding $x_jd_id = a_{ij}d$ to both sides, we have
$$x_i(d_jd - (-1)^d dd_j) + x_j (d_id - (-1)^d dd_i) = a_{ij} d,$$
or in other words
$$x_i[d_j,d] + x_j[d_i,d] = a_{ij}d.$$
It is easy to see that $d$ satisfies this equation if $d$ commutes with $d_i, d_j$ and $a_{ij}d = 0.$ Hence, 
$$Cent(d_1, \dots, d_n) \cap Ann(a_{ij}, i, j = 1, \dots, n) \subseteq Der(\Gamma_n(\mathcal{D})).$$


\end{document}